\documentclass[12pt]{article}
\usepackage{amsfonts}
\usepackage{epsfig}

\title{The Pentagram Integrals on Inscribed Polygons}
\author{Richard Evan Schwartz \thanks{\hskip 5 pt Supported by 
N.S.F. Research Grant DMS-0072607}\ \  and Serge Tabachnikov}

\newtheorem{theorem}{Theorem}[section]

\newtheorem{lemma}[theorem]{Lemma}

\def\startproof{{\bf {\medskip}{\noindent}Proof: }}

\def\endproof{$\spadesuit$  \newline}

\def\P{\mbox{\boldmath{$P$}}}%
\def\R{\mbox{\boldmath{$R$}}}%

\begin{document}
\maketitle

\section{Introduction} \label{intro}

The pentagram map is a geometric iteration defined on polygons. This
map is defined in practically any field, but it is most easily
described for polygons in the real projective plane.
Geometrically, the pentagram map carries the polygon $P$ to the
polygon $Q$, as shown in Figure 1. 

\begin{center}
\resizebox{!}{2.6in}{\includegraphics{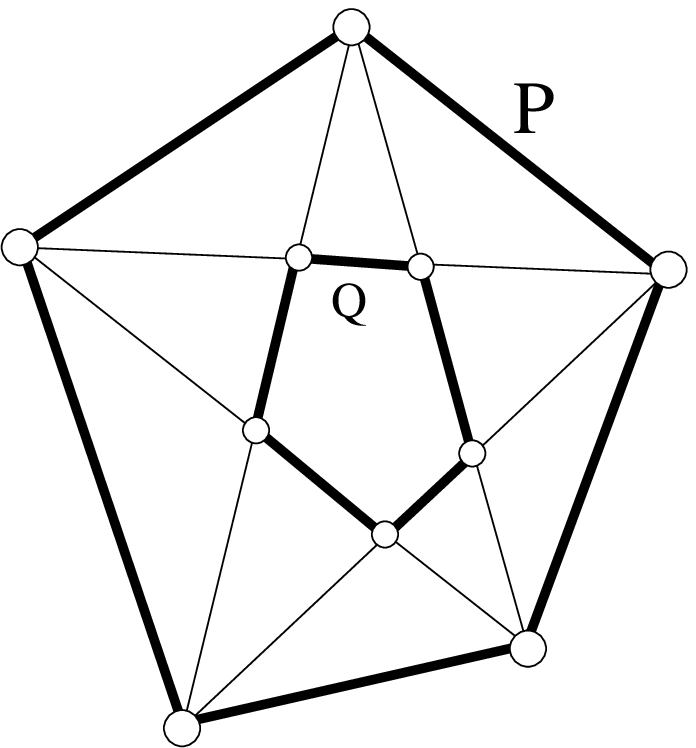}}
\newline
{\bf Figure 1:\/} The pentagram map
\end{center}

This pentagram map is always defined for convex polygons, and 
generically defined for all polygons. 
The pentagram map commutes with projective transformations and
induces a generically defined map $T$ on the space
${\cal Q\/}_n$ of  cyclically 
labelled  \footnote{Technically, one needs to consider the
square of the map in order to get a canonically defined
map on labelled $n$-gons.  However, if one is willing
to break symmety, preferring (say) left over right, then 
the map itself is defined on labelled $n$-gons.}
 $n$-gons modulo
projective transformations.   
$T$ is the identity map on ${\cal Q\/}_5$ and 
has period $2$ on ${\cal Q\/}_6$. The map
$T$ is not periodic on ${\cal Q\/}_n$
for $n \geq 7$.

The pentagram map also defines a map on a
slightly larger space ${\cal P\/}_n$.
The space ${\cal P\/}_n$ is the space of
{\it twisted $n$-gons\/} modulo projective
transformations.  A twisted $n$-gon is a
map $\phi: \R \to \R\P^2$, such that
$$
\phi(k+n)=M \circ \phi(k) \hskip 30 pt \forall k.
$$
Here $M$ is a real projective transformation,
called the {\it monodromy\/} of the twisted
polygon.  The space ${\cal P\/}_n$ contains
the space ${\cal Q\/}_n$ as a codimension $8$
submanifold.

This pentagram map is
studied for pentagons in [{\bf M\/}].  
We studied the general case in  [{\bf S1\/}],
[{\bf S2\/}], and [{\bf S3\/}].  Finally,
in [{\bf OST1\/}] (see [{\bf OST2\/}] for a short version) 
we showed that the pentagram map
 is a completely integrable system on
${\cal P\/}_n$.  We think that the pentagram
map is a completely integrable system on
${\cal Q\/}_n$ as well, but there are certain
technical details we need to resolve before
establishing this fact.  

In [{\bf S3\/}] we constructed polynomial functions
$$
O_1,...,O_{[n/2]},O_n,E_1,...,E_{[n/2]},E_n: {\cal P\/}_n \to \R
$$
which are invariant under the pentagram map.
In [{\bf OST1\/}] these functions are recognized a complete list of 
integrals for the completely integrable system.
There is an invariant Poisson structure, and the above
functions Poisson-commute with respect to this
Poisson structure.  We call these functions the
{\it monodromy invariants\/}, and we will construct
them in the next section.

The purpose of this paper is to prove the following result.
\begin{theorem}
\label{OE}
$O_k(P)=E_k(P)$ for any
twisted $n$-gon $P$ that is inscribed in a conic section and
any $k=1,...,[n/2], n$.
\end{theorem}

The result, in particular, applies to closed inscribed polygons.
Theorem \ref{OE} boils down to a countably infinite family
of polynomial identities.  The polynomials involved are
somewhat reminiscent of the symmetric polynomials,
but they have somewhat less symmetry then these.
We noticed this result by doing computer experiments. 
One novel feature of the theorem is that we discovered not
just the result but also the proof by way of computer
experimentation.   We wrote a Java applet
to aid us with the combinatorics of the proof.
This applet is available on the first author's 
website. \footnote{http://www.math.brown.edu/$\sim$res}

While our proof is mainly combinatorial, we think that
perhaps there should be a better proof based on geometry.
Accordingly, we will describe our polynomials
in three ways -- geometrically, combinatorially, and
in terms of determinants.  We will only use the
combinatorial description in our proof, but we
hope that the other descriptions might ring a
bell for some readers.

The pentagram map seems to interact nicely
with polygons that are inscribed in conic
sections.   We mention our paper
[{\bf ST\/}], in which we describe some
finite configuration theorems, \`a la Pappus,
that we discovered in connection with
the pentagram map and inscribed polygons.

 Here is an overview of the paper.
In \S 2 we define the monodromy invariants. In \S 3 we
reduce Theorem \ref{OE} to a combinatorial problem.
In \S 4 we solve this combinatorial problem. 

The second author would like to thank Brown University for its
hospitality during his sabbatical; he was also supported by NSF grant DMS-0072607.

\section{The Monodromy Invariants} \label{monoinv}

The papers [{\bf S3\/}] and [{\bf OST1\/}] give a good account
of the monodromy invariants.  The notation in the two
papers is slightly different, and we will follow [{\bf OST1\/}].
We will give three descriptions of these invariants, one
geometrical, one combinatorial, and one based on 
determinants.    The reader only
interested in the proof of Theorem \ref{OE} need only pay attention
to the combinatorial definition.  As we say in the introduction,
we mention the other definitions just in the hopes that it
will ring a bell for some reader.  

Because we are not
using the geometric or determinental definitions in
our proof, we will not include the arguments that
show the equivalence of the various definitions.
The paper [{\bf S3\/}] has a proof that the
geometric and combinatorial definitions coincide.

\subsection{The Geometric Definition} \label{geomdef}

First of all, we have the {\it cross ratio\/}
$$
[t_1,t_2,t_3,t_4]=\frac{(t_1-t_2)(t_3-t_4)}{(t_1-t_3)(t_2-t_4)}.
$$

Suppose that $\phi$ is a twisted $n$-gon, with monodromy $M$.  
We let   $v_i=\phi(i)$. The label $i$ in Figure 2 denotes $v_i$,
and similarly for the other labels.

\begin{center}
\resizebox{!}{2.8in}{\includegraphics{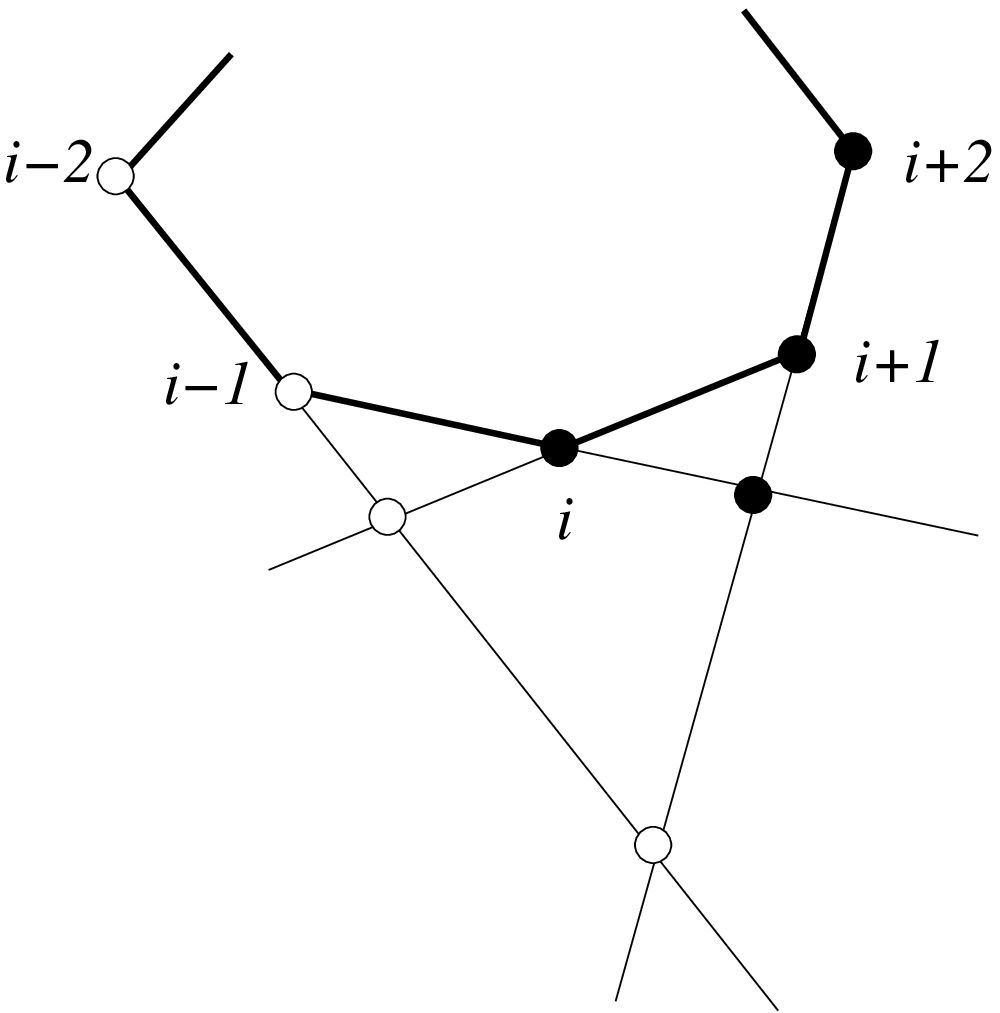}}
\newline
{\bf Figure 2:\/} vertex labels
\end{center}

We associate to each vertex $v_i$ two numbers:
$$
x_i=[v_{i-2},v_{i-1},((v_{i-2},v_{i-1}) \cap (v_i,v_{i+1}),((v_{i-2},v_{i-1}) \cap (v_{i+1},v_{i+2})],
$$
$$
y_i=[v_{i+2},v_{i+1},((v_{i+2},v_{i+1}) \cap (v_i,v_{i-1}),((v_{i+2},v_{i+1}) \cap (v_{i-1},v_{i-2})].
$$
Here $(a,b)$ denotes the line determined by points $a$ and $b$.
For instance, $x_i$ is the cross ratio of the $4$ white points in Figure 2.
We call the invariants $x_1,y_1,x_2,y_2,...$ the {\it corner invariants\/}.
These invariants form a periodic sequence of length $2n$. That is
$x_{k+n}=x_k$ and $y_{k+n}=y_k$ for all $k$.

We define
$$
O_n=\prod_{i=1}^n x_i; \hskip 30 pt
E_n=\prod_{i=1}^n y_i.
$$

The other monodromy invariants are best defined in an indirect way.
Recall that $M$ is the monodromy of our twisted polygon $\phi$.
We lift $M$ to an element of $GL_3(\R)$ which we also denote by
$M$.  We define
$$
\Omega_1=\frac{{\rm trace\/}^3(M)}{{\rm det\/}(M)}; \hskip 30 pt
\Omega_2=\frac{{\rm trace\/}^3(M^{-1})}{{\rm det\/}(M^{-1})}.
$$
These quantities are independent of the lift of $M$ and only
depend on the conjugacy class of $M$.  Finally, these quantities
are invariant under the pentagram map.

We define
$$
\widetilde \Omega_1=O_n^2E_n \Omega_1; \hskip 30 pt
\widetilde \Omega_2=O_nE_n^2 \Omega_2.
$$
It turns out that these quantities are polynomials in
the corner invariants. The remaining monodromy invariants
are suitably weighted homogeneous parts of these polynomials.

We have a basic rescaling operation
$$
R_t(x_1,y_1,x_2,y_2,...)=(tx_1,t^{-1}y_1,tx_2,t^{-1}y_2,...).
$$
We say that a polynomial in the corner invariants has {\it weight $k$\/}
if
$$
R_t^*(P)=t^kP.
$$
Here $R_t^*$ denotes the obvious action of $R_t$ on polynomials.
In [{\bf S3\/}] we show that
$$
\widetilde \Omega_1=\sum_{k=1}^{[n/2]} O_k ; \hskip 30 pt
\widetilde \Omega_2=\sum_{k=1}^{[n/2]} E_k.
$$
Here $O_k$ and $E_k$ are the weight $k$ polynomials
in each sum, and $[n/2]$ denotes the floor of $n/2$.

\subsection{The Combinatorial Definition} 
\label{combdef}

Now we describe the combinatorial formulas for these invariants.
Again, the paper [{\bf S3\/}] has a proof that the
description below coincides with the one in the previous section.
In everything we say, the indices are taken cyclically, mod $n$.

We introduce the monomials 
\begin{equation}
\label{Xx}
X_i:=x_i\,y_i\,x_{i+1}.
\end{equation}

The monodromy polynomial $O_k$ is built out
of the monomials $x_i$ for $i=1,...,n$ and
$X_j$ for $j=1,..,n$.   
We call two monomials {\it consecutive\/} if they
 involve consecutive or coinciding variables $x_i$.
More explicitly, we have the following criteria:
\begin{enumerate}
\item $X_i$ and $X_j$ are consecutive if 
$
j\in\left\{i-2,\,i-1,\,i,\,i+1,\,i+2\right\};
$
\item
$X_i$ and $x_j$ are consecutive if
$j\in\left\{i-1,\,i,\,i+1,\,i+2\right\};$
\item $x_i$ and $x_j$ are consecutive if $j \in \{i-1,i,i+1\}$.
\end{enumerate}

Let $O$ be a monomial obtained by the product of the monomials
$X_i$ and $x_j$, that is, 
$$
 O=X_{i_1}\cdots{}X_{i_s}\,x_{j_1}\cdots{}x_{j_t}. 
$$
Such a monomial is called {\it admissible} if no two of the monomials
are consecutive.
For every admissible monomial, define the {\it weight} $|O|=s+t$
and the {\it sign} $\mathrm{sign}(O)=(-1)^s$.
One then has
$$
O_k=\sum_{|O|=k}\mathrm{sign}(O)\,O; \hskip 30 pt 
k\in\left\{1,2,\ldots,\left[\frac{n}{2}\right]\right\}.
$$

For example, if $n \geq 5$ we obtain the following two polynomials:
$$
O_1=\sum_{i=1}^n
\left(
x_i-x_i\,y_i\,x_{i+1}
\right),
\qquad
O_2=\sum_{i=1}^n
\left(
x_i\,x_{i+2}-
x_i\,y_i\,x_{i+1}\,x_{i+3}
\right),
$$

The same formula works for $E_k$, if we make all the same
definitions with $x$ and $y$ interchanged. More precisely, one builds the polynomials $E_k$ from the monomials $y_i$ and $Y_j:=y_{j-1}x_jy_j$ with the same restriction that no consecutive monomials are allowed and the same definitions of the weight and sign.

We note that the dihedral symmetry 
\begin{equation}
\label{dih1}
\sigma(x_i)=y_{-i},\ \sigma(y_i)=x_{-i}
\end{equation}
interchanges the polynomials $O_k$ and $E_k$.

\subsection{The Determinantal Definition}
\label{detdef}

Now we describe determinantal formulas for the monodromy invariants; these formulas did not appear in our previous papers on the subject.  

For positive integers $k>l$, we define the four-diagonal determinant
$$
F_l^k=\left|\begin{array}{ccccccc}
1&x_k&X_{k-1}&0&0&\dots&0\\
1&1&x_{k-1}&X_{k-2}&0&\dots&0\\
0&1&1&x_{k-2}&X_{k-2}&\dots&0\\
\dots&\dots&\dots&\dots&\dots&\dots&\dots\\
0&\dots&\dots&\dots&1&1&x_{l+1}\\
0&\dots&\dots&\dots&0&1&1
\end{array}\right|,
$$
where $x_i$ and $y_i$ are the corner invariants and $X_i$ is as in (\ref{Xx}). By convention,
$$
F_{k+2}^k=0, \ F_{k+1}^k=1, \ F_{k}^k=1.
$$
Then one has the following formula for the monodromy invariants $O_k$:
\begin{equation} \label{allO}
\sum_{i=0}^{[n/2]} O_i=F_1^n+F_0^{n-1}-F_1^{n-1}+x_ny_nx_1F_2^{n-1}.
\end{equation}

Similarly, for $E_k$, define 
$$
G_p^q=\left|\begin{array}{ccccccc}
1&y_{p+1}&Y_{p+2}&0&0&\dots&0\\
1&1&y_{p+2}&Y_{p+3}&0&\dots&0\\
0&1&1&y_{p+3}&Y_{p+4}&\dots&0\\
\dots&\dots&\dots&\dots&\dots&\dots&\dots\\
0&\dots&\dots&\dots&1&1&y_{q}\\
0&\dots&\dots&\dots&0&1&1
\end{array}\right|.
$$
Then one has:
\begin{equation} \label{allE}
\sum_{i=0}^{[n/2]} E_i= G_0^{n-1}+G_1^n-G_1^{n-1}+y_nx_1y_1G_1^{n-2}.
\end{equation}

Formulas (\ref{allO}) and (\ref{allE}) simplify if one considers an {\it open} version of the monodromy invariants: instead of having a periodic ``boundary condition" $x_{i+n}=x_i, y_{i+n}=y_i$, assume that $x_i=0$ for $i\leq 0$ and $i\geq n+1$, and $y_i=0$ for $i\leq -1$ and $i\geq n$. With this ``vanishing at infinity" boundary conditions, the monodromy invariants are given by a single determinant:
$$
\sum_{i=0}^{[n/2]} O_i=F_0^n,\quad \sum_{i=0}^{[n/2]} E_i=G_{-1}^{n-1}.
$$ 
We shall encounter an open version of monodromy invariants in Section \ref{open}.
\section {Reduction to the Puzzle} \label{reduct}

A non-degenerate conic in $\R\P^2$ can be identified with $\R\P^1$ by 
way of the stereographic projection from a point of the conic. This 
identification is unique  up to a projective transformation of the real 
projective line. That is, a different choice of the center of projection 
amounts to a projective transformation of  $\R\P^1$. 

Recall from the introduction that
a twisted $n$-gon is a
map $\phi: \R \to \R\P^2$, such that
$$
\phi(k+n)=M \circ \phi(k) \hskip 30 pt \forall k.
$$
Here $M$ is a real projective transformation,
called the {\it monodromy\/} of the twisted
polygon. 
Given a twisted $n$-gon, we let $v_j=\phi(j)$.
These are the vertices of the twisted $n$-gon.

If $(\dots,v_{-1},v_0,v_1,\dots)$ is an inscribed twisted polygon, 
we can consider the vertices $v_i$ as points of the real projective 
line and talk about their cross-ratios which are uniquely defined. 
Referring to the cross-ratio on the projective line, 
we set: 
$$
p_i=1-[v_{i-2},v_{i-1},v_i,v_{i+1}].
$$
In the next lemma we express the corner invariants 
of an inscribed polygon in terms of the quantities $p_i$.
Once we specify the cross ratios $\{p_i\}$, we produce
an inscribed twisted polygon having corner
invariants $\{(x_i,y_i)\}$.  Thus we have a map
$(p_i)\mapsto (x_i,y_i)$.  We denote this map by $F$.

\begin{lemma} \label{xycross}
One has:
$$
x_i=[v_{i-2},v_{i-1},v_i,v_{i+2}],\quad y_i=[v_{i-2},v_i,v_{i+1},v_{i+2}]
$$
and the map $F$ is given by the formula
\begin{equation} 
\label{xypeq}
x_i=\frac{1-p_{i}}{p_{i+1}},\quad y_i=\frac{1-p_{i+1}}{p_{i}}.
\end{equation}
\end{lemma}

\startproof
Consider Figure 3. Using the projection from point $v_{i+1}$, we obtain:
$$
x_i=[v_{i-2},v_{i-1},A,B]=[(v_{i+1}v_{i-2}),(v_{i+1}v_{i-1}),(v_{i+1}A),(v_{i+1}B)]=[v_{i-2},v_{i-1},v_i,v_{i+2}].
$$
A similar projection from $v_{i-1}$ yields the formula for $y_i$.
The expression for $x_i$  in terms of $p_i$ and $p_{i+1}$ follows now from the identity
$$
[v_{i-2},v_{i-1},v_i,v_{i+2}]=\frac{[v_{i-2},v_{i-1},v_i,v_{i+1}]}{1-[v_{i-1},v_{i},v_{i+1},v_{i+2}]},
$$
and likewise for $y_i$.
\endproof

\begin{center}
\resizebox{!}{4.8in}{\includegraphics{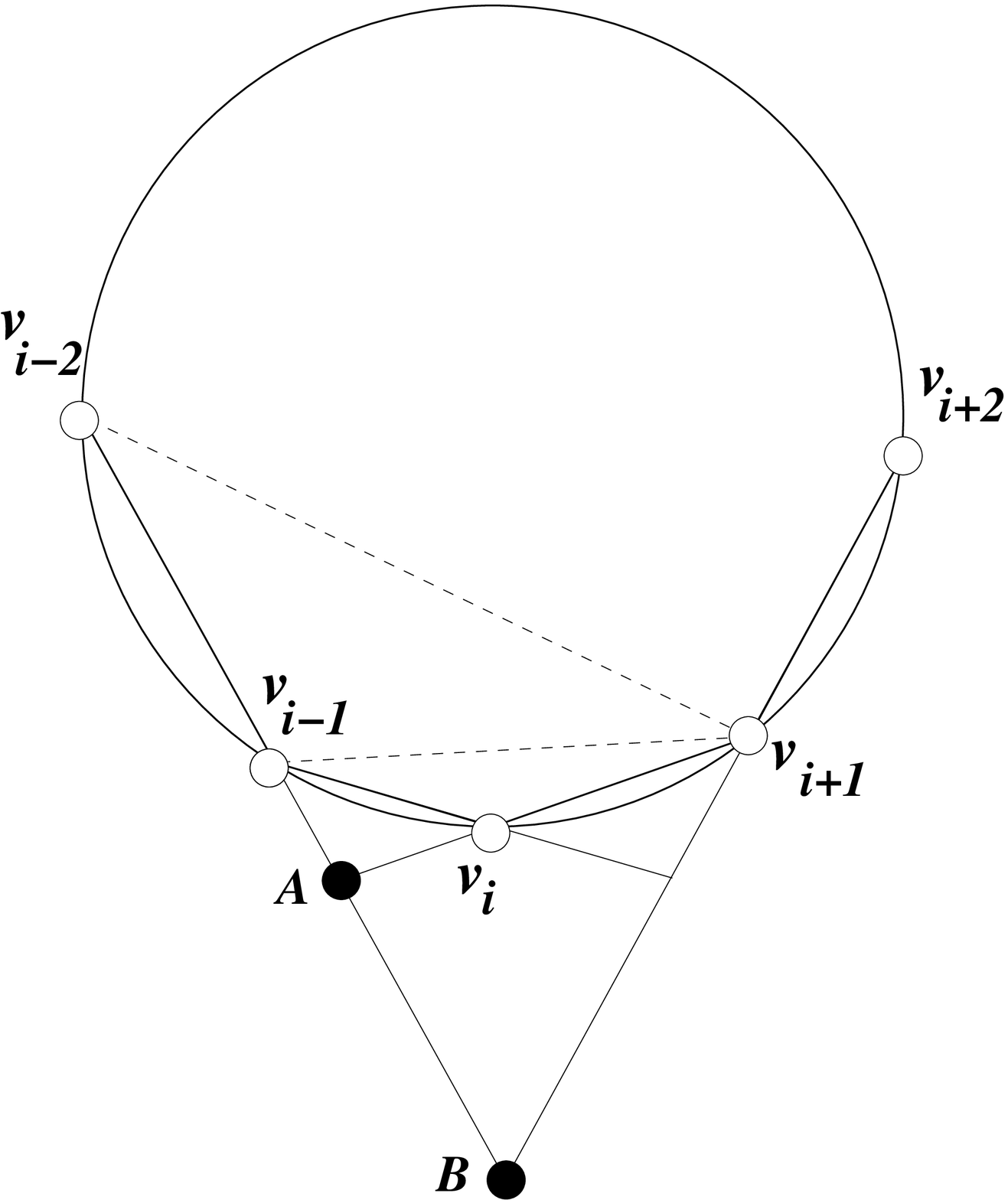}}
\newline
{\bf Figure 3:\/} Proof of Lemma \ref{xycross}
\end{center}

The dihedral group also acts on cyclic sequences $(p_i)$ by cyclic 
permutations and the orientation-reversing involution  $\sigma'(p_i)=p_{1-i}$. 
It follows from (\ref{xypeq}) that  $F\circ \sigma'=\sigma\circ F$ where 
$\sigma$ is the involution (\ref{dih1}). Hence $F$ is a dihedrally equivariant map.

After $x_i$ and $y_i$ are replaced by $p_i$ via (\ref{xypeq}), the 
polynomials $O_k$ and $E_k$ become Laurent polynomials in the 
variables $p_i$. The identity $E_n=O_n$ obviously holds since 
both sides equal $\prod (1-p_i)/p_i$. We need to prove that 
$E_k=O_k$ for $k=1,...,[n/2]$. The strategy is to show that 
every monomial in the variables $p_i$ appears in $E_k$ and 
$O_k$ with the same coefficient. Since the map
$F$ is dihedrally equivariant 
and the orientation-reversing involution on the variables $(x_i,y_i)$ 
interchanges $E_k$ and $O_k$, it suffices to show that if two monomials 
in the variables $p_i$ are related by orientation-reversing symmetry 
(for example, the involution $\sigma'$) then they appear in $O_k$ with the same coefficients. 

Let us compute the monomials in the variables $p$'s that appear
 in the polynomial $O_k$. Using (\ref{xypeq}), one finds:
\begin{equation}
\label{xp}
x_i=\frac{1}{p_{i+1}}-\frac{p_i}{p_{i+1}}
\end{equation}
and
\begin{equation}
\label{Xp} 
-X_i=-\frac{1}{p_ip_{i+1}p_{i+2}}+
\frac{1}{p_{i+1}p_{i+2}}+
2\frac{1}{p_{i}p_{i+2}}-2\frac{1}{p_{i+2}}
-\frac{p_{i+1}}{p_{i}p_{i+2}}
+\frac{p_{i+1}}{p_{i+2}}.
\end{equation}

We see that the variables $p$'s that appear in $x_i$ involve two 
indices, $i$ and $i+1$, and the those  in $X_i$ involve three 
indices, $i,i+1$ and $i+2$. Pictorially, these terms can be 
represented as follows: for $x_i$, see Figure 4, and for $X_i$, see Figure 5.
In these figures, the presence of each term $p$ in the numerator 
and denominator is represented by a shaded square, and its absence by a white square. 

\begin{center}
\resizebox{!}{1.2in}{\includegraphics{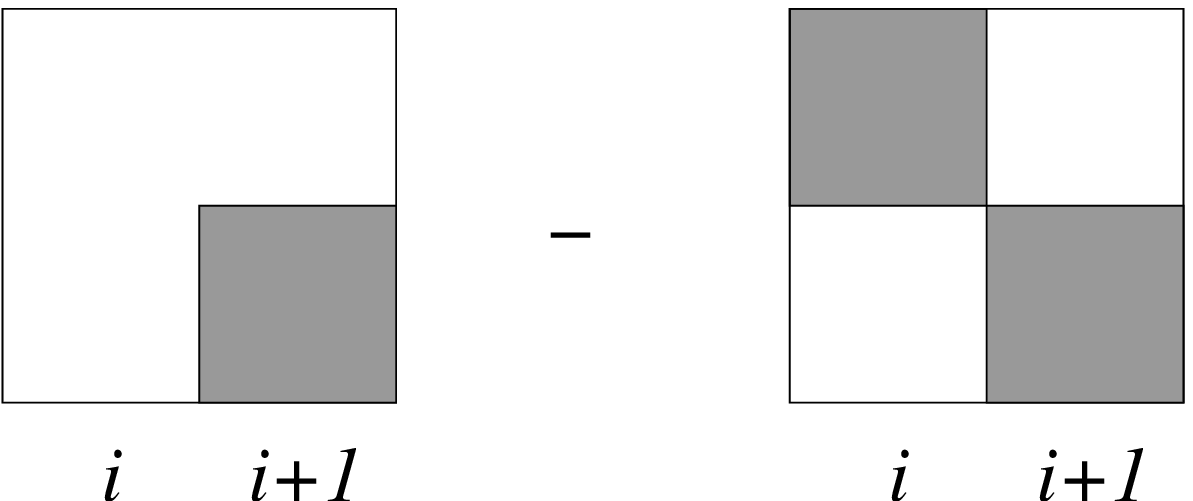}}
\newline
{\bf Figure 4:\/} pictorial representation of (\ref{xp})
\end{center}

\begin{center}
\resizebox{!}{2.4in}{\includegraphics{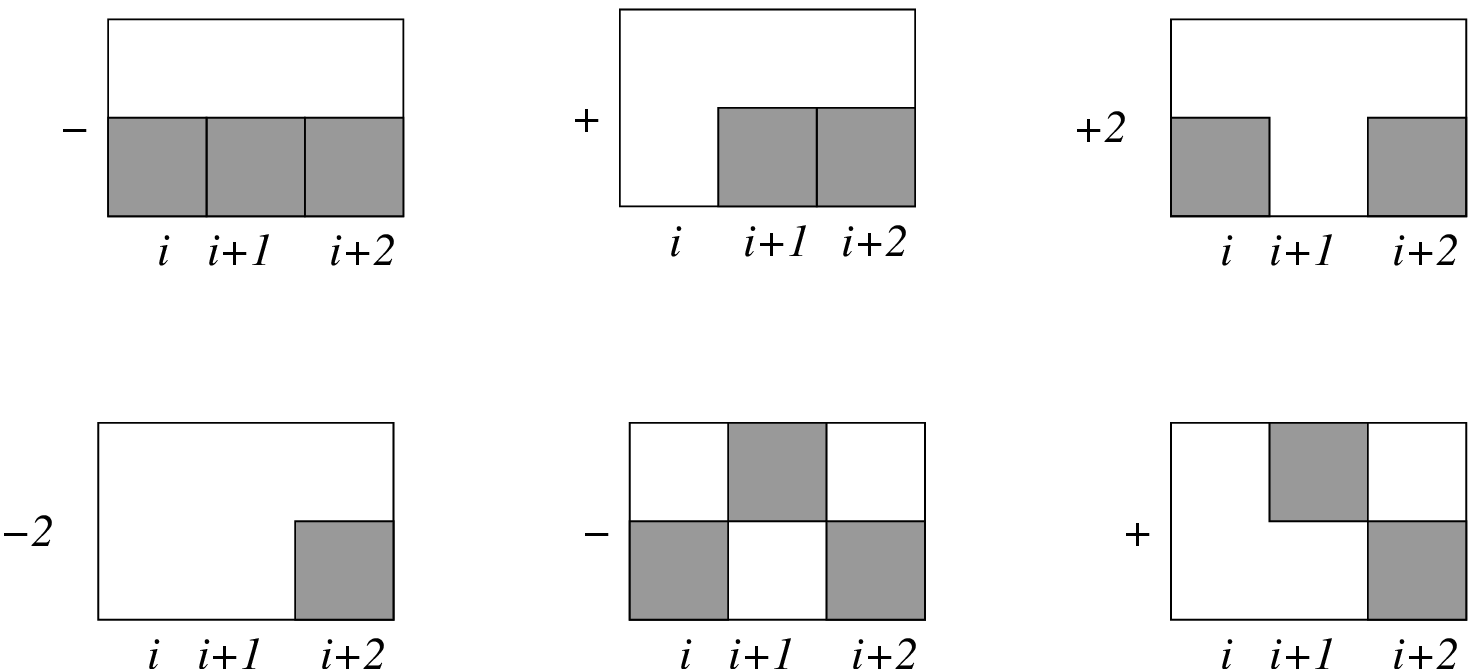}}
\newline
{\bf Figure 5:\/} pictorial representation of (\ref{Xp})
\end{center}

According to Section \ref{combdef}, the polynomial $O_k$ is the 
sum of all admissible products of $k$ terms, and each term is 
either $x_i$ or $(-X_j)$. The admissibility condition is that the 
monomials $x_i$ or $X_j$ involved are sufficiently distant; 
what it means precisely is that the respective tiles in Figures 4 and 5, 
corresponding to these terms, do not overlap. This is a crucial observation.

To summarize, each monomial in $O_k$, after the 
substitutions (\ref{xp}) and (\ref{Xp}), is represented by a 
collection of $k$ tiles depicted in Figures 4 and 5, taken with 
the product of their respective coefficients. The tiles, that occupy 
two or three consecutive positions,  are placed around a circle having 
$N$ positions available (if we are concerned with twisted $N$-gons). 
There may be empty positions left between the tiles.

As a final preparation for the next section,  we introduce the following 
notation: letter $A$ denotes a shaded square in the lower row of a tile ($p_i$ in denominator), 
letter $B$ a shaded square in the upper row ($p_i$ in numerator), and letter $X$ an empty 
column in a tile. Thus, the two tiles in Figure 4 correspond to the 
words $XA$ and $BA$, and the six tiles in Figure 5 
to $AAA, XAA, AXA, XXA, ABA$ and $XBA$, respectively. 
We also use letter $X$ to mean an empty slot between the tiles. 

In the next section we reformulate Theorem \ref{OE} 
as a statement about
a kind of puzzle involving words in letters $\{A,X,B\}$.
After we formulate the combinatorial statement, we prove it.
The combinatorial result we prove implies
Theorem \ref{OE}.

\section{The Puzzle}

\subsection{The Main Result}

Now we are going to extract the main combinatorial
information from the discussion at the end of the
last section.

We fix some integer $N>0$ and consider the set of
length $N$ lists in the letters $\{A,X,B\}$.    We consider
two lists the same if they are cyclic permutations of each other.
We say that a {\it cyclic sentence\/} is
an equivalence class of such strings.    To illustrate our notation by
way of example, 
$(AXABA)$ denotes the equivalence class of $AXABA$.  Here $N=5$.
We let $\cal S$ denote the set of such sentences, with the
dependence on $N$ understood.

We single out certain strings of letters, and to each of these special
{\it words\/} we assign a coefficient and a weight.
Here is the list.
$$
\matrix{{\rm word\/} & {\rm coefficient\/} & {\rm weight\/} \cr
X & 1 & 0 \cr
XA & 1 & 1 \cr
XAA & 1 & 1 \cr
XBA & 1 & 1 \cr
AXA & 2 & 1 \cr
AAA & -1 & 1 \cr
BA & -1 & 1 \cr
ABA & -1 &1 \cr
XXA & -2 & 1 \cr
}
$$

We say that a {\it parsing\/} of a cyclic sentence is a description of the
cyclic sentence as a concatenation of words.  For example 
$$(ABA/XA); \hskip 30pt (XAA/BA)$$
are the only two parsings of $(ABAXA)$.  
We define the {\it coefficient\/} of a parsing to be
the product of the coefficients of the word.  We define
the {\it weight\/} of the parsing to be the
sum of the weights of the words.   Both parsings
above have coefficient $-1$ and weight $2$.

For each cyclic sentence $S$, we define
$c(S,w)$ to be the total sum of the coefficients
of the weight $w$ parsings of $S$.  For instance,
when $S=(ABAXA)$, we have
$c(S,2)=-2$ and otherwise $c(S,w)=0$.
For a more streamlined equation, we define
$$
|S|=\sum_{w=0}^{\infty} c(S,w)t^w.
$$
Here $|S|$ is a polynomial in $t$ that encodes all the
coefficients at once.  For example
$$|(ABAXA)|=-2t^2.$$

Let $\overline S$ denote the cyclic sentence obtained by reversing $S$.
In view of Section \ref{reduct}, 
Theorem \ref{OE} is equivalent to the following result.
\begin{theorem}[Cyclic]
\label{cyclic}
We have $|S|=|\overline S|$ for all
cyclic sentences $S$.
\end{theorem}

In light of the work in the previous section,
Theorem \ref{cyclic} implies Theorem \ref{OE}.
The rest of this chapter is devoted to proving
Theorem \ref{cyclic}.

\subsection{The Tight Puzzle}

Before we tackle Theorem \ref{cyclic}, we slightly 
modify our puzzle for the sake of convenience.

\begin{lemma}
\label{trivial}
Suppose that $(W)$ contains the string $XB$.  Then
$|W|=0$.
\end{lemma}

\startproof
No word ends in $B$, and so the string $XB$ must
continue as $XBA$.
The occurence of $.../X/BA/...$ in any parsing contributes
weight $1$ and coefficient $-1$ whereas the occurence of
$.../XBA/...$ contributes weight $1$ and coefficient $1$.
If we have a  parsing that involves
$.../XBA/...$ we can create a new parsing by replacing
the last $.../XBA/...$ with $.../X/BA/...$.  These two
parsings have the same weight and opposite coefficient,
and thereby cancel each other in the total sum.
\endproof

By Lemma \ref{trivial}, we can
simply throw out any strings that contain
$XB$, and we may drop the word
$XBA$ from our list of words.    
There is a similar cancellation involving $XXA$.
Within a parsing, the occurence of
$.../X/XA/...$ has weight $1$ and coefficient $1$ whereas the
occurence of $.../XXA/...$ has weight $1$ and
coefficient $-2$.
If we have a  parsing that involves
$.../XXA/...$ we can create a new parsing by replacing
the last $.../XXA/...$ with $.../X/XA/...$.   The new parsing
cancels out ``half'' of the original.   Thus, we may consider
an alternate puzzle where the parsing $/X/XA/$ is forbidden and
the word list is

$$
\matrix{{\rm word\/} & {\rm coefficient\/} & {\rm weight\/} \cr
X & 1 & 0 \cr
XA & 1 & 1 \cr
XAA & 1 & 1 \cr
AXA & 2 & 1 \cr
AAA & -1 & 1 \cr
BA & -1 & 1 \cr
ABA & -1 &1 \cr
XXA & -1 & 1 \cr
}
$$
All we have done is dropped $XBA$ from the list and
changed the weight of $XXA$ from $-2$ to $-1$.
We call this last puzzle the {\it tight puzzle\/}.   
Establishing Theorem \ref{cyclic}
for the tight puzzle is the same as establishing these results
for the original one.

\subsection{The Open Version}
\label{open0}

As an intermediate step to proving Theorem \ref{cyclic}, we state a
variant of the result.
We consider bi-infinite strings
in the letters $\{A,B,X\}$, where there are only finitely many
$A$s and $B$s.     We say that two such strings are {\it equivalent\/}
if one of them is a shift of the other one.  We say that an
{\it open sentence\/} is an equivalence class of such strings.
We use finite strings to denote sentences, with the understanding
that the left and right of the finite string is to be padded with
an infinite number of $X$s.  Thus, $ABAXA$ refers to the
bi-infinite sentence
$...XXABAXAXX...$.  We define parsings just as in the cyclic case.
For instance, here are all the parsings of this sentence
\begin{itemize}
\item $/ABA/XA/(-1)$
\item $/XXA/BA/XA/(1)$
\end{itemize}
(Recall that we have forbidden $/X/XA$.)
The first  of these have weight $2$ and the last one has weight $3$.
We have put the coefficients next to the parsings in each case.
Our notation is such that the left and right sides of each expression
are padded with $.../X/X/...$.     Based on the list above, we have
$$|ABAXA|=-t^2+t^3.$$

Here is the variant of Theorem \ref{cyclic} for open sentences.

\begin{theorem}[Open]
\label{open}
We have $|S|=|\overline S|$ for all
open sentences $S$.
\end{theorem}

Theorem \ref{cyclic} implies Theorem \ref{open} in a fairly direct way.
For instance, suppose we are interested in proving Theorem \ref{open}
for an open sentence $S$.    We say that the {\it span\/} of $S$ is
the combinatorial distance between the first and last non-$X$ letter
of $S$.  For instance, the span of $ABAXA$ is $4$.   Supposing that
$S$ has span $s$, we simply create a cyclic sentence of length
(say) $s+10$ by padding the nontrivial part of $S$ with $X$s
and then taking the cyclic equivalence class.  Call this cyclic
sentence $S'$.   We clearly have
$$|S|=|S'|=|\overline S'|=|\overline S|.$$
The middle equality is Theorem \ref{cyclic}.
The end inequalities are obvious.

Now we turn to the proof of Theorem \ref{open}.
In the next result, $W$ stands for a finite string in the
letters $A,B,X$.

\begin{lemma}[Right Identities]
\label{right}
The following identities hold.
\begin{enumerate}
\item $|WAAA|+t|W|=0$.
\item $|WXAA|-t|W|=0$.
\item $|WXXA|+t|W|=0$.
\item $|WAXA|-t|WA|-2t|W|=0$.
\item $|WABA|+t|WA|+t|W|=0$.
\end{enumerate}
\end{lemma}

\startproof
Consider Identity 1.
Any parsing of $WAAA$ must have the form
$W/AAA$.    But $AAA$ has weight $1$ and
coefficient $-1$.   Hence
$c(WAAA,w)=-c(W,w-1)$.   Also
$c(WAAA,0)=0$.  Identity 1 follows immediately from this.
Identity 2 and Identity 3 have the same proof.

Consider Identity 4.
There are two kinds of parsings of $WAXA$.  One kind
has the form $W/AXA$ and the other kind has the form
$/WA/XA$.    Note that $AXA$ has weight $1$ and
coefficient $2$ and $XA$ has weight $1$ and
coefficient $1$.   From this, we see that
$c(WAXA,w)=2c(W,w-1)+c(WA,w-1)$ for
all $w$.  Identity 4 follows immediately.
Identity 5 has the same proof.
\endproof

\noindent
{\bf Discussion:\/}
If Theorem \ref{open} really holds, then the
``reverses'' of all the identities above should always hold.
Let's consider an example in detail.
The reverse of Identity 2 above is
$$|AAXW|-t|W|=0,$$
for all strings $W$.   However, taking $W=ABA$, the
weight $3$ parsings of $AAXW$ are
\begin{itemize}
\item $/XXA/AXA/BA/ (2)$.
\item $/XAA/XA/BA/ (-1)$.
\end{itemize}
As usual, our convention is to leave off the words $.../X/X/...$ on both
sides.    Adding up the coefficients, we see that
$c(AAXW,3)=1$.     At the same time, 
the only weight $2$ parsing of $W$ is
\begin{itemize}
\item $/XXA/BA/ (1)$
\end{itemize}
Hence $c(W,2)=1$.     This accords with our supposed equality,
but the $3$ parsings in the one case don't obviously cancel
out the $2$ parsings in the other.   
In Lemma \ref{right},
the various parsings matched up and cancelled each other
in an obvious way.  However, this does not happen
for the reverse identities.    Nonetheless, we will 
prove the reverse identities of Lemma \ref{right}.

\begin{lemma}[Left Identities]
\label{left}
The following identities hold.
\begin{enumerate}
\item $|AAAW|+t|W|=0$.
\item $|AAXW|-t|W|=0$.
\item $|AXXW|+t|W|=0$.
\item $|AXAW|-t|AW|-2t|W|=0$.
\item $|ABAW|+t|AW|+t|W|=0$.
\end{enumerate}
\end{lemma}

\startproof
We will prove Identity 5.  The other identities have the same proof.
First of all, we check computationally that Identity 5 holds (say)
for all strings $W$ having length at most $3$.

Suppose now that $W$ is a shortest word for which we don't
know the result of this lemma.  We know that $W$ has length
at least $3$, so we can write $W=VR$, where $R$ has length
$3$ and ends in $A$.  Consider the case when $R=AXA$.

By induction, we have
\begin{equation}\label{eq1}
|ABAV|+t|AV|+t|V|=0.\end{equation}
\begin{equation}\label{eq2}
|ABAVA|+t|AVA|+t|VA|=0.\end{equation}
Using Identity 4 of Lemma \ref{right} (three times) we have
\begin{equation}
\label{eq3}
t|VAXA|=t^2|VA|+2t^2|V|.
\end{equation}
\begin{equation}
\label{eq4}
t|AVAXA|=t^2|AVA|+2t^2|AV|.
\end{equation}
\begin{equation}
\label{eq5}
|ABAVAXA|=t|ABAVA|+2t|ABAV|.
\end{equation}
When we add together the right hand sides of Equations \ref{eq3}, \ref{eq4}, \ref{eq5},
we get $0$, thanks to Equations \ref{eq1} and \ref{eq2}.
Hence, when we add the left hand sides of Equations \ref{eq3}, \ref{eq4}, \ref{eq5},
we also get $0$.  But this last sum is exactly the identity we wanted to prove.

A similar argument works when $R$ is any of the other $3$-letter strings
that appear in Lemma \ref{right}.     The only case we haven't considered
is the case when $R=XBA$, but these strings are forbidden.
\endproof

Now that we have Lemma \ref{right} and Lemma \ref{left}, our proof of
Theorem \ref{open} goes by induction.  
First of all, we check Theorem \ref{open} for all strings having span
at most $3$.   Suppose then that $W$ is the
shortest open sentence   for which we do not know Theorem \ref{open}.
We can write $W=VR$ where $R$ is some string of length $3$ that
ends in $A$.    

Let's consider the case when $R=XAA$.   Then we have
$$|W|=|VXAA|=t|V|=t|\overline V|=|AAX\overline V|=|\overline W|.$$
The second equality is Identity 2 of Lemma \ref{right}.  The third equality is the 
induction assumption.  The fourth equality is Identity 2 of Lemma \ref{left}.
A similar argument works when $R$ is any of the $3$ letter strings
in Lemma \ref{right}.    The final case, $R=XBA$, is forbidden.

This completes the proof of Theorem \ref{open}.

\subsection{The Cyclic Case}

We need to mention another convention before we launch into
the cyclic case.
Besides cyclic and open sentences, there is one more case
we can consider.  We introduce the notation $[W]$ to denote
an open word whose parsings cannot be created by padding $X$s onto
the left and right of $W$.
We will illustrate what we have in mind
by way of example.     Setting $W=ABAXA$,  the parsings of
the open string $W$ are
$$/AXA/BA/; \hskip 30 pt
/XXA/X/ABA.$$
However, the second parsing involves two $X$s that have
been padded onto the left of $W$.  Only the first
parsing of $W$ is also a parsing of the {\it locked} string $[W]$.
We let $|[W]|$ be the polynomial that encodes the
weights and coefficients of all the parsings of $[W]$.

\begin{lemma}
Theorem \ref{cyclic} holds for any cyclic word $W$ with
no $X$ in it.
\end{lemma}

\startproof
To avoid some messy notation, we will consider an example.  The example
is sufficiently complex that it should illustrate the general proof. 
Suppose that
\begin{equation}
\label{exp}
W=(BA^2BA^5BA^1BA^7).
\end{equation}
Here, for instance $A^2=AA$.
Any parsing of $W$ must have the breaks
$$W=(BA/A^1BA/A^4BA/A^0BA/A^6).$$
The point is that we must have a break after each $BA$.
From this, we see that
$$|W|=|[A^6BA]| \times |[ABA]| \times |[A^4BA]| \times |[A^0BA]|.$$
To get the list of exponents on the right hand side of this product,
we simply decrement each exponent in Equation \ref{exp} by one.
But we would get the same list of exponents (perhaps in a different
order) when considering the reverse word $\overline W$.
\endproof

Below we prove the following result.
\begin{lemma}
\label{mainID}
The relation
$$|(WX)|+2|(WB)|-|W|=0$$
holds for all open words $W$.
\end{lemma}

Lemma \ref{mainID} allows us to finish the proof of Theorem \ref{cyclic}.   Our proof goes
by induction on the number of $Xs$ in the word $W$.   Let
$W$ be a word having the smallest number of $X$s, for which
we do not know Theorem \ref{cyclic}.     

After cyclically permuting the letters in $W$, we can
write $W=VX$. By Lemma \ref{mainID}, we have
$$|(W)|=|(VX)|=-2|(VB)|+|V|.$$
Applying Lemma \ref{mainID} to $\overline W$, we have
$$|(\overline W)=|(\overline VX)|=-2|(\overline VB)|+|\overline V|=
-2|\overline{VB}|+|\overline V|.$$
Setting $Y=VB$, we have
$$
|(W)|=-2|(Y)|+|V|= -2|(\overline Y)|+|\overline V| = |(\overline W)|.
$$
The middle equality comes from Theorem \ref{open} (to handle $V$) and the
induction assumption (to handle $Y$.)

\subsection{Some Auxilliary Relations} \label{aux}

It remains only to prove Lemma \ref{mainID}.  We will establish some
auxilliary relations in this section, and then use them 
in the next section to prove Lemma \ref{mainID}.

\begin{lemma} \label{taut}
For any string $W$, 
$$|XXW|=|W|.$$
\end{lemma}

\startproof
This is a tautology.
\endproof

\begin{lemma} \label{useful}
For any string $W$, 
$$
|W|=|[XW]|-|[AXW]|.
$$
\end{lemma}

\startproof
The proof makes use of the right identities from Lemma \ref{right} and is 
similar to the proof of Theorem \ref{open}. One first checks the 
statement  for all strings of span three, and then argues inductively 
on the span. The induction step is proved using the right identities 
from Lemma \ref{right} that hold verbatim for locked strings as well.

To illustrate the idea, we assume that $W$ ends with $AXA$, that is, 
$W=VAXA$. Then, by Identity 4 of Lemma \ref{right}, 
$$
|W|=|VAXA|=t|VA|+2t|V|,
$$
$$
|[XW|=|[XVAXA]|=t|[XVA]|+2t|[XV]|,
$$
and
$$
|[AXW]|=|[AXVAXA]|=t|AXVA|+2t|[AXV]|.
$$
By the induction assumption, 
$$
|VA|=|[XVA]|-|[AXVA]|,\ \ |V|=|[XV]|-|[AXV]|,
$$
and the result follows for $W$.
\endproof

\begin{lemma}
\label{axa}
For any string $W$, 
$$
|(AXAW)|=2t|[W]|+|[XAWA]|.
$$
\end{lemma}

\startproof
Either a parsing of the cyclic word $(AXAW)$ contains the 
string $AXA$, or there is a break after the first $A$ in  $(AXAW)$. The 
former case corresponds to the first term, $2t|[W]|$, 
and the latter case to the second, $|[XAWA]|$.
\endproof

\begin{lemma}
\label{aba}
For any string $W$,
$$
|(ABAW)|=|[WABA]|.
$$
\end{lemma}

\startproof
There must be a break after second $A$ in $(ABAW)$, and this
 provides a one-to-one correspondence between the parsings of $(ABAW)$ and $[WABA]$.
\endproof

\subsection{Proof of Lemma \ref{mainID}} \label{finaltouch}

\begin{lemma}
Lemma \ref{mainID} holds if $W$ does not start and end with $A$.
\end{lemma}

\startproof
This is a case-by-case analysis.
Suppose that
 $W=VX$ for some word $V$. Then $(WX)=(VXX)$. Since we must 
have a break between $V$ and $XX$, it follows that $$|(VXX)|=|XXV|=|V|=|W|.$$
(the second equality holds by Lemma \ref{taut}). On the other hand, $$|(WB)|=|(VXB)|=0$$ 
since the combination $X/BA$ is prohibited. The claim of Lemma \ref{mainID} follows.
Similarly, if $W=XV$ for some word $V$ then $(WX)=(VXX)$, and the same argument applies. 
If $W=VB$ then then each term in the equality of 
Lemma \ref{mainID} vanishes. The same holds if $W=BV$.
Finally, if $W$ ends in $B$, all terms in Lemma \ref{mainID} are trivial.
\endproof

The only remaining case is when $W=AVA$ for some word $V$. What we need to prove is 
\begin{equation}
\label{final}
|(AXAV)|+2|(ABAV)|-|AVA|=0.
\end{equation}

By Lemma \ref{axa}, 
$$|(AXAV)|=2t|[V]|+|[XAVA]|.$$
By Lemma \ref{aba}, 
$$|(ABAV)|=|[VABA]|.$$
 By Lemma \ref{useful}, 
$$|AVA|=|[XAVA]|-|[AXAVA]|.$$
 Therefore, the left hand side of (\ref{final}) equals
\begin{equation}
\label{veryfinal}
2t|[V]| + |[XAVA]|+ 2|[VABA]| - |[XAVA]| + |[AXAVA]|.
\end{equation}
By Identity 5 of Lemma \ref{right} for locked words, 
$$|[VABA]|=-t|[VA]|-t|[V]|.$$
 Finally, 
$$|[AXAVA]|=2t|[VA]|$$
 since 
a parsing of $[AXAVA]$ must start with $AXA$.  It follows that (\ref{veryfinal}) equals
$$
2t|[V]| + |[XAVA]| -2t|[VA]|-2t|[V]|- |[XAVA]| + 2t|[VA]|=0,
$$
as needed.  

This completes the proof of Lemma \ref{mainID}.

\section{References}

\noindent
[{\bf M\/}] Th. Motskin, {\it The Pentagon in the projective plane, and a
comment on Napier's rule\/}, Bull. Amer. Math. Soc. {\bf 52\/} (1945), 985-989.
\newline
\newline
[{\bf OST1\/}] V. Ovsienko, R. E. Schwartz, S. Tabachnikov,
{\it The Pentagram Map: A Discrete Integrable System\/},
Comm. Math. Phys. 2010 (to appear),  
math arXiv:0810.5605. 
\newline
\newline
[{\bf OST2\/}] V. Ovsienko, R. E. Schwartz, S. Tabachnikov,
{\it Quasiperiodic motion for the pentagram map\/},
Electr. Res. Announ. Math. 
{\bf 16} (2009), 1--8.
\newline
\newline
[{\bf S1\/}] R. Schwartz, {\it The Pentagram Map\/} 
J. Experiment. Math. {\bf 1\/} (1992),  71--81.
\newline
\newline
[{\bf S2\/}] R. Schwartz, {\it Recurrence of the Pentagram Map\/},
J. Experiment. Math. {\bf 110\/} (2001),  519--528.
\newline
\newline
[{\bf S3\/}] R. Schwartz, {\it Discrete Monodromy, Pentagrams,
 and the Method of Condensation\/}, 
J. Fixed Point Theory   Appl. {\bf 3\/} (2008),  379--409.
\newline
\newline
[{\bf ST\/}] R. Schwartz, S. Tabachnikov,
{\it Elementary Surprises in Projective Geometry\/}, 
Math. Intelligencer, 2010 (to appear), math arXiv:0910.1952.

\end{document}